\documentclass[journal,twoside]{IEEEtran}

\usepackage{fancyhdr}
\usepackage{cite}
\usepackage[dvips]{graphicx}
\pdfoutput=1
\usepackage{amsmath,amssymb,amsfonts}
\usepackage{algorithmic}
\usepackage{array}
\usepackage[caption=false,font=footnotesize]{subfig}
\usepackage{stfloats}
\usepackage{url}
\usepackage{textcomp}
\usepackage{xcolor}
\usepackage{lipsum}
\def\BibTeX{{\rm B\kern-.05em{\sc i\kern-.025em b}\kern-.08em
    T\kern-.1667em\lower.7ex\hbox{E}\kern-.125emX}}
\usepackage{adjustbox}
\usepackage{multirow}
\usepackage{hhline}
\usepackage{tablefootnote}
\usepackage{booktabs} 
\usepackage{color}
\usepackage{tikztensor}
\definecolor{kulblue}{RGB}{0,85,165}
\tikzset{fancy/.style={inner color=kulblue!5!white, outer color=kulblue!20!white, fill faces, opacity=0.95}}
\colorlet{kulblue20}{kulblue!20!white}
\colorlet{kulblue70}{kulblue!70!white}
\colorlet{kulblue30}{kulblue!30!black}
\colorlet{kulblue60}{kulblue!60!black}
\colorlet{kulblue90}{kulblue!90!black}

\newcommand{\bc}[1]{\mbox{\boldmath $\mathcal{#1}$}}
\newcommand{\mb}[1]{\mathbf{#1}}
\newcommand{\F}{\mathrm{F}}
\newcommand{\T}{\mathrm{T}}

\begin{document}

\title{Block-Term Tensor Decomposition: Model Selection and Computation}

\author{Athanasios~A.~Rontogiannis,~\IEEEmembership{Member,~IEEE,}
        Eleftherios~Kofidis,~\IEEEmembership{Member,~IEEE,}
        and~Paris~V.~Giampouras%
\thanks{A. A. Rontogiannis is with the IAASARS, National Observatory of Athens, 152~36 Penteli, Greece. E-mail: tronto@noa.gr.}%
\thanks{E. Kofidis is with the Dept. of Statistics and Insurance Science, University of Piraeus, 185~34 Piraeus, Greece and the Computer Technology Institute \& Press ``Diophantus" (CTI), 265~04 Patras, Greece. E-mail: kofidis@unipi.gr}%
\thanks{P. V. Giampouras is with the Mathematical Institute for Data Science, Johns Hopkins University, Baltimore, MD 212~18, USA. E-mail: 
parisg@jhu.edu}
}

\markboth{}%
{Rontogiannis \MakeLowercase{\textit{et al.}}: Block-Term Tensor Decomposition: Model Selection and Computation}

\maketitle

\begin{abstract}
The so-called block-term decomposition (BTD) tensor model has been recently receiving increasing attention due to its enhanced ability of representing systems and signals that are composed of \emph{blocks} of rank higher than one, a scenario encountered in numerous and diverse applications. Its uniqueness and approximation have thus been thoroughly studied. Nevertheless, the challenging problem of estimating the BTD model structure, namely the number of block terms and their individual ranks, has only recently started to attract significant attention. In this paper, a novel method of BTD model selection and computation is proposed, based on the idea of imposing column sparsity \emph{jointly} on the factors and in a \emph{hierarchical} manner and estimating the ranks as the numbers of factor columns of non-negligible magnitude. Following a block successive upper bound minimization (BSUM) approach for the proposed optimization problem is shown to result in an alternating hierarchical iteratively reweighted least squares (HIRLS) algorithm, which is fast converging and enjoys high computational efficiency, as it relies in its iterations on small-sized sub-problems with closed-form solutions. Simulation results for both synthetic examples and a hyper-spectral image de-noising application are reported, which demonstrate the superiority of the proposed scheme over the state-of-the-art in terms of success rate in rank estimation as well as computation time and rate of convergence.
\end{abstract}


\IEEEpeerreviewmaketitle

\section{Introduction}
\label{sec:intro}

\IEEEPARstart{B}{lock}-Term Decomposition (BTD) was introduced in~\cite{ldl08b} as a tensor model that combines the Canonical Polyadic Decomposition (CPD) and the Tucker decomposition (TD), in the sense that it decomposes a tensor in a sum of tensors that have low multilinear rank
(instead of rank one as in CPD\footnote{Note that a rank-1 tensor is also a rank-$(1,1,\ldots,1)$ tensor.}). In other words, BTD is a sum of TDs (block terms). Hence a BTD can be seen as a \emph{constrained} TD, with its core tensor being block diagonal (see \cite[Fig.~2.3]{ldl08b}). Given the sum-of-TDs structure of BTD and in view of the fact that CPD is also a constrained TD~\cite{sdfhpf17}, BTD can also be seen as a \emph{constrained} CPD having factors with (some) collinear columns~\cite{ldl08b}. In a way, BTD lies between the two extremes (in terms of core tensor structure), CPD and TD, and it is interesting to recall the related remark made in~\cite{ldl08b}, namely that ````the" rank of a higher-order tensor is actually a combination of the two aspects: one should specify the number of blocks \emph{and} their size." 
Accurately and efficiently estimating these numbers for a given tensor is the main subject of this work.

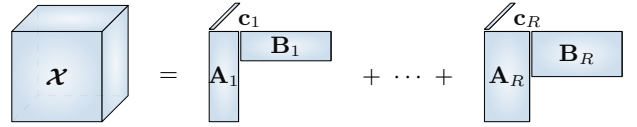
\begin{figure}
\centering
\begin{tikzpicture}[node distance=0.3cm, chain,
term/.style={dim={4,4,3},fancy,tensor scale=0.3}]
\small
\node [tensor, term, dashed back lines] {$\bc{X}$};
\node {$=$};
\node [rank-LL1 tensor, term, L=1.3, 2d,
labels={$\mathbf{A}_{1}$}{$\mathbf{B}_{1}$}{$\mathbf{c}_{1}$}] {};
\node {$+\ \cdots\ +$};
\node [rank-LL1 tensor, term, L=2, 2d,
labels={$\mathbf{A}_{R}$}{$\mathbf{B}_{R}$}{$\mathbf{c}_{R}$}] {};
\end{tikzpicture}
\caption{Rank-$(L_r,L_r,1)$ block-term decomposition.}
\label{fig:BTD}
\end{figure}

Although~\cite{ldl08b} introduced BTD as a sum of $R$ rank-$(L_r,M_r,N_r)$ terms ($r=1,2,\ldots,R$) in general, the special case of rank-$(L_r,L_r,1)$ BTD has attracted a lot more of attention, because of both its more frequent occurrence in applications and the existence of more concrete and easier to check uniqueness conditions. This paper will also focus on this special yet very popular BTD model. Consider a 3rd-order tensor, $\bc{X}\in\mathbb{C}^{I\times J\times K}$. Then its rank-$(L_r,L_r,1)$ decomposition is written as
\begin{equation}
\bc{X}=\sum_{r=1}^{R}\mb{E}_{r}\circ \mb{c}_{r},
\label{eq:BTD}
\end{equation}
where $\mb{E}_{r}$ is an $I\times J$ matrix of rank $L_r$, $\mb{c}_{r}$ is a nonzero column $K$-vector and $\circ$ denotes outer product. 
Clearly, $\mb{E}_{r}$ can be written as a matrix product $\mb{A}_{r}\mb{B}_{r}^{\T}$ with the matrices $\mb{A}_{r}\in\mathbb{C}^{I\times L_r}$ and $\mb{B}_{r}\in\mathbb{C}^{J\times L_r}$ being of full column rank, $L_r$. A schematic diagram of the rank-$(L_r,L_r,1)$ BTD is shown in Fig.~\ref{fig:BTD}.

BTD has found applications in communications (e.g., \cite{lb08,sel17}), neuro- and anatomical imaging~\cite{hcspvhl14,ckmt19,ahmrv19,scsh19}, electrocardiography (ECG)~(e.g., \cite{z17,oz19,msf19}), hyper-spectral imaging (HSI)~\cite{qxzzt17,xzq19,zfhw19,bg20}, community detection in networks~\cite{gpp20}, spectrum cartography~\cite{zfwzh20}, and electron microscopy~\cite{srp17}, among others. Recently it has also been proposed as a compact model of neural networks in modern machine learning applications~\cite{mzzdhsz19,ylcyzx20}. 
The application of BTD in blind source separation (BSS) was first considered in~\cite{ldl11} and later presented in more detail in~\cite{ldl12}, giving rise to the so-called Block Component Analysis (BCA) approach. The underlying idea is that BTD can better represent components (sources) of a variable complexity (hence rank), while CPD-based BSS\footnote{also referred to as Canonical Polyadic Analysis (CPA)~\cite{ldl12}.} restricts the sources to have rank one.\footnote{An intuitively pleasant way to describe this difference is to say that, while CPA decomposes the data into ``atoms", BCA decomposes it into ``molecules"~\cite{ldl12}.}

The uniqueness of BTD was studied in~\cite{ldl08b}, also for the general rank-$(L_r,M_r,N_r)$ case. Essential uniqueness for the rank-$(L_r,L_r,1)$ BTD of eq.~(\ref{eq:BTD}) means that the only indeterminacies are the order of the $R$ terms and a scaling of the $\mb{E}_{r}$ matrix with a counter-scaling of the vector $\mb{c}_{r}$.
The most popular (though not the only one) uniqueness theorem for this case states that a \emph{sufficient} uniqueness condition is that the partitioned matrices $\mb{A}\triangleq \left[\begin{array}{cccc} \mb{A}_1 & \mb{A}_2 & \cdots & \mb{A}_R \end{array}\right]$ and $\mb{B}\triangleq \left[\begin{array}{cccc} \mb{B}_1 & \mb{B}_2 & \cdots & \mb{B}_R \end{array}\right]$ are of full column rank and $\mb{C}\triangleq \left[\begin{array}{cccc} \mb{c}_1 & \mb{c}_2 & \cdots & \mb{c}_R \end{array}\right]$ does not have any collinear columns~\cite[Theorem~4.1]{ldl08b}. The generic version of the requirement for full column rank of $\mb{A},\mb{B}$ is that $\min(I,J)\geq \sum_{r=1}^{R}L_r$, which can easily be met in applications where $R$ and $L_r$ are small. It should however be noted that this is not a necessary condition as our simulation results also demonstrate. 

Alternating least squares (ALS) was 
extended to the computation of a tensor BTD in~\cite{ldl08c}. In that same work, it was also shown (and demonstrated through an example) that degeneracy can also occur for BTD.\footnote{That a best BTD approximation of given ranks may not exist for a real-valued tensor was later shown in~\cite{gc19}. This is not the case for tensors in $\mathbb{C}$, however (cf.~\cite{gc19} and references therein).} In the noise-free case, and as shown in~\cite[Theorem~4.1]{ldl08b}, the BTD can be also computed with the aid of a generalized eigenvalue decomposition (GEVD), provided the above uniqueness condition is satisfied. Recently, algebraic solution methods that are free from this limitation have been also reported~\cite{dl19,dl20,vdd20}. In the presence of noise, these solutions can serve to initialize the ALS iterations~\cite{ldl08c}. ALS with the appropriate modifications to incorporate the non-negativity constraint was used in~\cite{qxzzt17} for non-negative BTD of hyper-spectral imagery. Non-alternating (all-at-once) computation approaches, including gradient descent and nonlinear least squares, were followed in~\cite{svd13} and the resulting methods are implemented in Tensorlab~\cite{vdsvd16}. Additional methods of BTD computation include ALS regularized through $\ell_2$ norms of its factors (to avoid over-fitting~\cite{zfwzh20} or to enforce low rank~\cite{xzq19}) or through proximal point modifications~\cite{Li2013}, deflation-based~\cite{gc17b}, variable projection using Riemannian gradient for rank-$(L_r,M_r,N_r)$ BTD with factors of orthonormal columns~\cite{oal18}, tensor block diagonalization~\cite{tpc17},  solving the equivalent matrix factorization problem with one of the factors constrained to have low-rank rows~\cite{fh19}, and computing an appropriately constrained coupled CPD~\cite{vdd20}.

In most of the BTD methods mentioned above, $R$ and $L_r$, $r=1,2,\ldots,R$ are assumed known (and it is commonly assumed that all $L_r$ are all equal to $L$, for simplicity). In fact, in practice, this is a challenging question on its own. Unless external information is given (such as in a telecommunications ~\cite{ldl12} or in a HSI unmixing application with given or estimated ground truth~\cite{qxzzt17}), there is no way to know these values a priori. An observation that is common in all known BCA applications is that the separation performance does not strongly depend on the particular values of the $L_r$ ranks~\cite{qxzzt17,ckmt19,zfwzh20}. In fact, as it was also observed in~\cite{hcspvhl14,ckmt19}, the method is robust to overestimation of $L_r$ (although, as observed in~\cite{scsh19}, performance of BTD-based classifiers may considerably vary with $L_r$). Nonetheless, one should try not to set $L_r$ to a very high value. The reason is that, in addition to increasing the computational complexity, setting $L_r$ too high may hinder interpretation of the results through letting noise/artifact sources interfere with the desired sources~\cite{hcspvhl14}. 
This holds for $R$ as well, although its choice is known to be more crucial to the obtained performance. For example, setting $R$ too high in~\cite{hcspvhl14} results in source splitting (also referred to as over-factoring~\cite{hakss17}), thus compromising the separation and interpretation of the components.

\subsection{Background}
\label{subsec:back}

Model order selection techniques for BTD can be dictated by corresponding CPD techniques, as reviewed in~\cite[Section~4]{hcspvhl14}, including clustering similar CPD components (e.g., \cite{tpc17}). Schemes of multilinear rank estimation (largely based on matrix rank estimation and/or extensions of one-dimensional information-theoretic criteria) are also relevant in view of the constrained TD structure of BTD~\cite{p16,ylc17,xzmz18,hakss18,ahmrv19}. 
In the absence of noise, the model rank parameters can be computed as a by-product of recently reported algebraic BTD methods~\cite{dl19,dl20}. Thus, in the non-iterative method of~\cite{dl20}, and in the (almost) noise-free case, these are estimated (with the aid of singular value decompositions (SVDs)) from a joint block matrix diagonalization problem. For noisier tensors, $R$ and $\sum_r L_r$ are assumed known.

Model order selection can also be application-specific. For example, $L_r$'s are estimated in~\cite{z17} as the auto-regressive (AR) orders of the sources in ECG analysis, with $R$ assumed known. In~\cite{ckmt19}, and in the functional magnetic resonance imaging (fMRI) context, $L_r$ is estimated as the number of statistically significant (bearing useful information) columns of $\mb{A}_r,\mb{B}_r$. \cite{xzq19} relies on the subspace-based method of~\cite{bn08} for estimating the number $R$ of spectral signatures in BTD-based HSI denoising.   

Alternative techniques rely on sparsity arguments for model selection. A greedy scheme, inspired from a sparse coding viewpoint, is proposed in~\cite{bppc13}, for more general tensor decompositions~\cite{pctzl13} including BTD as a special case. Instead of building the model incrementally, however, one can follow the reverse way of starting from a rank \emph{overestimate} and arrive at the true rank(s) by eliminating negligible components, aided in this task by appropriate regularization. Such an approach is followed in~\cite{hakss17,h19}, where the constrained CPD formulation of BTD is taken advantage of to first estimate $R$ and then $L_r$'s assumed all equal, before computing the model factors in~(\ref{eq:BTD}). In each case, a regularization term is added to the tensor approximation cost, which is composed of mixed norms of the factor matrices and serves as upper bound on the tensor nuclear norm thus promoting column sparsity to the factors and hence low rank. The augmented Lagrangian method is adopted for the computations. 

Nevertheless, as demonstrated in~\cite{yws16},\cite{tgrk18} for the CPD case, the problems of model rank estimation and approximation of factors can be addressed \emph{jointly}, with significant gains in both accuracy and complexity (of particular interest for big data applications). This idea is proposed in~\cite{gofzc20} for the rank-$(L_r,L_r,1)$ BTD model with not necessarily all equal block-term ranks $L_r$. A regularization term consisting of the sum of the mixed $\ell_{1,2}$ norms of the matrices $\mb{A},\mb{B},\mb{C}$ is added to the squared error of the tensor approximation, namely
\begin{align}
&\underset{\mb{A},\mb{B},\mb{C}}{\min}\frac{1}{2}\left\|\bc{Y} - \sum_{r=1}^R\mb{A}_r\mb{B}_r^{\T}\circ\mb{c}_r\right\|_{\mathrm{F}}^2 + \nonumber \\
 &\gamma (\|\mb{A}\|_{1,2}+\|\mb{B}\|_{1,2}+\|\mb{C}\|_{1,2}),
 \label{eq:ABC}
 \end{align}
where $\|\cdot\|_{\F}$ is the Frobenius norm,  $\|\cdot\|_{1,2}$ denotes the mixed $\ell_{1,2}$ norm (defined as the $\ell_1$ norm of the $\ell_2$ norms of the matrix columns\footnote{In~\cite{gofzc20}, this is referred to as the $\ell_{2,1}$ norm.}), and $\gamma$ is the regularization parameter weighing the regularization term over the data fidelity term. This sparsity-inducing regularization helps promoting low rank for the BTD factors and hence estimating $R$ (as the number of non-zero columns of $\mb{C}$) and $L_r$'s (as the number of non-zero columns of the $r$th blocks of $\mb{A},\mb{B}$ that correspond to non-zero columns of $\mb{C}$). For the solution of~(\ref{eq:ABC}), a proximal term is first added in~\cite{gofzc20} and then a block coordinate descent (BCD) approach is taken, leading to a regularized version of the ALS procedure of~\cite{ldl08c} that will be referred to henceforth as the BTD alternating group lasso (BTD-AGL) algorithm.

The approach we propose in this paper also falls in the previous category. Yet, it has a number of very important new features, inherited from our earlier work on factorization-based low-rank approximation of matrices~\cite{grk16},\cite{grk17},\cite{grk19a}, from which it draws inspiration.
In~\cite{grk19a}, the sum of reweighted Frobenius norms of the factors of the data matrix is used as regularization and, in particular, a diagonal weighting, jointly depending on the factors, is proposed, naturally leading to an iteratively reweighted least squares (IRLS)~\cite{ddfg09} solution approach, with fast convergence and low complexity.
Here we generalize that idea in the BTD problem. The regularization of~\cite{grk19a} is employed, in two levels: first, combining the reweighted norms of $\mb{A}$ and $\mb{B}$, and second, coupling these with the reweighted norm of $\mb{C}$. This two-level coupling naturally matches the structure of the model in~(\ref{eq:BTD}), making explicit the different roles of $\mb{A},\mb{B}$ and $\mb{C}$, in contrast to previous related works~\cite{hakss17,gofzc20} that miss to exploit this relation. Furthermore, due to this fact, the regularization proposed here has a stronger sparsity promoting action compared with previous works. Applying majorization with appropriate upper bounds and a BCD approach results in an alternating hierarchical IRLS (HIRLS) algorithm that manages to both reveal the ranks and compute the BTD factors at a high convergence rate and low computational cost. Notably, iterations involve closed-form updates that contain only matrix-matrix 
multiplications, which can be efficiently implemented on most modern computer systems and are easily parallelizable. The complexity can be reduced even more by eliminating negligible columns (column pruning) in the course of the iterations (as in~\cite{grk19a} for the low-rank matrix factorization problem). Simulation results for both synthetic examples and a HSI de-noising application are reported, which demonstrate the superiority of the proposed scheme over the state-of-the-art in terms of success rate in rank estimation as well as computation time and rate of convergence.\footnote{A short version of this work was accepted for presentation in EUSIPCO-2020.} 

The rest of this paper is organized as follows. The adopted notation is described in the following subsection. The problem is mathematically stated in Section~\ref{sec:problem}. The proposed method is presented in Section~\ref{sec:method}. Section~\ref{sec:examples} reports and discusses the simulation results. Conclusions are drawn and future work plans are outlined in Section~\ref{sec:concls}.

\subsection{Notation}
\label{subsec:notation}

Lower- and upper-case bold letters are used to denote vectors and matrices, respectively. Higher-order tensors are denoted by upper-case bold calligraphic letters. For a tensor $\bc{X}$, $\mb{X}_{(n)}$ stands for its mode-$n$ unfolding. $\otimes$ stands for the Kronecker product. The Khatri-Rao product is denoted by $\odot$ in its general (partition-wise) version and by $\odot_{\mathrm{c}}$ in its column-wise version. $\circ$ denotes the outer product. The superscript $^{\T}$ stands for transposition. The identity matrix of order $N$ and the all ones column $N$-vector are respectively denoted by $\mb{I}_N$ and $\mb{1}_N$. The row vectorization and the trace of a matrix $\mb{X}$ are denoted by $\mathrm{vec}(\mb{X})$ and $\mathrm{tr}(\mb{X})$, respectively. $\nabla_{\mb{X}}$ stands for the gradient operator with respect to (w.r.t) $\mb{X}$. $\mathrm{diag}(\mb{x})$ is the diagonal matrix with the vector $\mb{x}$ on its main diagonal. The Euclidean vector norm and the Frobenius matrix norms are denoted by $\|\cdot\|_{2}$ and $\|\cdot\|_{\F}$, respectively. The mixed 1, 2 ($\ell_{1,2}$) norm of a matrix $\mb{X}=\left[\begin{array}{ccc} \mb{x}_1 & \cdots & \mb{x}_N \end{array}\right]$ is defined as $\sum_{n=1}^{N}\|\mb{x}_i\|_2$. $\mathbb{C}$ is the field of complex numbers. 

\section{Problem Statement}
\label{sec:problem}

Given an $I\times J\times K$ tensor $\bc{Y}$, its best (in the least squares sense) rank-$(L_r,L_r,1)$ approximation is sought for, namely
\begin{equation}
\underset{\mb{A},\mb{B},\mb{C}}{\min}f(\mb{A},\mb{B},\mb{C})\triangleq \frac{1}{2}\left\|\bc{Y} - \sum_{r=1}^R\mb{A}_r\mb{B}_r^{\T}\circ\mb{c}_r\right\|_{\F}^2,
\label{eq:BTDapprox}
\end{equation}
where the matrices $\mb{A}_r=\left[\begin{array}{cccc} \mb{a}_{r1} & \mb{a}_{r2} & \cdots & \mb{a}_{rL_r} \end{array}\right]\in\mathbb{C}^{I\times L_r}$, 
$\mb{B}_r=\left[\begin{array}{cccc} \mb{b}_{r1} & \mb{b}_{r2} & \cdots & \mb{b}_{rL_r} \end{array}\right]\in\mathbb{C}^{J\times L_r}$, $\mb{C}\in\mathbb{C}^{K\times R}$, and the ranks $R$ and $L_r$, $r=1,2,\ldots,R$ are unknown. In terms of its mode unfoldings, the tensor $\bc{X}\triangleq\sum_{r=1}^R\mb{A}_r\mb{B}_r^{\T}\circ\mb{c}_r$ can be written as~\cite{ldl08b}
\begin{eqnarray}
\mb{X}_{(1)}^{\T} \!\!\!\!\! & = & \!\!\!\!\! (\mb{B}\odot\mb{C})\mb{A}^{\T}, \label{eq:X1} \\
\mb{X}_{(2)}^{\T} \!\!\!\!\! & = & \!\!\!\!\! (\mb{C}\odot\mb{A})\mb{B}^{\T}, \label{eq:X2} \\
\mb{X}_{(3)}^{\T} \!\!\!\!\!\!\! & = & \!\!\!\!\!\!\!
\left[\begin{array}{ccc} (\mb{A}_1\odot_{\mathrm{c}} \mb{B}_1)\mb{1}_{L_1} & \cdots & (\mb{A}_R\odot_{\mathrm{c}} \mb{B}_R)\mb{1}_{L_R}\end{array}\right]\mb{C}^{\T}.
\label{eq:X3}
\end{eqnarray}
These expressions can be used in alternatingly solving for $\mb{A},\mb{B},\mb{C}$, respectively.

The regularization-based approach adds terms to the objective function above with the aim of imposing constraints on the sought factors, as in~(\ref{eq:ABC}) for example. However, in contrast to~(\ref{eq:ABC}), where all BTD factors are treated in the same manner, the regularizer proposed in this paper perfectly matches the structure of the BTD model, offering increased flexibility via a suitable joint block and column sparsity promoting mechanism of a hierarchical nature. The proposed modification to~(\ref{eq:BTDapprox}) can be stated as
\begin{equation}
\underset{\mb{A},\mb{B},\mb{C}}{\min}f(\mb{A},\mb{B},\mb{C})+\lambda \|\mb{F}(\mb{A},\mb{B},\mb{C})\|_{1,2},
\label{eq:regularized}
\end{equation}
where regularization is performed with the aid of the $\ell_{1,2}$ norm of the $2\times R$ matrix $\mb{F}(\mb{A},\mb{B},\mb{C})$, constructed as follows. Let $\mb{G} \triangleq \left[\begin{array}{cc} \mb{A}^\T &  \mb{B}^\T\end{array}\right]^\T$ be the $(I+J)\times\sum_{r=1}^{R}L_r$ matrix resulting from stacking the factors $\mb{A}$ and $\mb{B}$ and $\mb{G}_r \triangleq \left[\begin{array}{cc} \mb{A}_r^\T&  \mb{B}_r^\T\end{array}\right]^\T$ denote its $r$th $(I+J)\times L_r$ block. The matrix $\mb{F}(\mb{A},\mb{B},\mb{C})$ is defined as
\begin{equation}
\mb{F}(\mb{A},\mb{B},\mb{C}) \triangleq \left[\begin{array}{cccc}
\|\mb{G}_1\|_{1,2} & \|\mb{G}_2\|_{1,2} & \cdots & \|\mb{G}_R\|_{1,2}  \\
\|\mb{c}_1\|_2 & \|\mb{c}_2\|_2 & \cdots &  \|\mb{c}_R\|_2
\end{array}\right].
\label{eq:matF}    
\end{equation}
The minimization of the $\ell_{1,2}$ norm of a vector or matrix subject to a data proximity criterion has been widely utilized in the literature for enforcing group sparsity in vector/matrix recovery problems~\cite{hlmqy17}. This property of the $\ell_{1,2}$ norm was exploited in our earlier work~\cite{grk19a,grk19b} for model order selection in low-rank matrix factorization applications. In the present work, we extend that idea to the BTD problem by employing a two-level hierarchical $\ell_{1,2}$ norm-based regularization scheme. At the upper level, the $\ell_{1,2}$ norm of the matrix $\mb{F}(\mb{A},\mb{B},\mb{C})$ above promotes the elimination of whole blocks of $\mb{A}$ and $\mb{B}$ (which are tied together by the mixed norms $\|\mb{G}_r\|_{1,2}$, $r=1,2,\ldots,R$) and the corresponding columns of $\mb{C}$. At the lower level, the $\ell_{1,2}$ norms $\|\mb{G}_r\|_{1,2}$ induce column sparsity to the ``surviving" blocks of $\mb{A}, \mb{B}$. Hence, we have the flexibility to  overestimate the ranks $R$ and  $L_r$, $r=1,2,\ldots,R$ as $R=R_{\mathrm{ini}}$ and $L_r=L_{\mathrm{ini}}$ in the unknown BTD model, since this regularization can reduce them towards their actual values with a proper selection of the regularization parameter $\lambda$. The problem in~(\ref{eq:regularized}) may be solved using a block successive upper bound minimization (BSUM) approach~\cite{hrlp16}, as described in the next section. As explained in Appendix~\ref{sec:B}, the resulting algorithm is an alternating hierarchical IRLS scheme, referred to henceforth as \emph{BTD-HIRLS}. 

\section{Proposed Method}
\label{sec:method}

First, we rewrite the minimization problem~(\ref{eq:regularized}) more explicitly in terms of the BTD factors $\mb{A},\mb{B}$, and $\mb{C}$ as \begin{align}
&\underset{\mb{A},\mb{B},\mb{C}}{\mathrm{min}}\frac{1}{2}\left\|\bc{Y} - \sum_{r=1}^R\mb{A}_r\mb{B}_r^{\T}\circ\mb{c}_r\right\|_{\F}^2 + \nonumber \\
 &\lambda\sum_{r=1}^{R}\sqrt{\left(\sum_{l=1}^{L}\sqrt{\|\mb{a}_{rl}\|_{2}^2 + \|\mb{b}_{rl}\|_{2}^2 + \eta^2}\right)^2 + \|\mb{c}_r\|_{2}^2 + \eta^2},
\label{eq:minp} 
\end{align}
where $\eta^2$ is a very small positive constant that ensures smoothness and $R$ and $L$ here stand for the initial {\it (over)estimates} of the model rank parameters. It can be shown that the objective function in~(\ref{eq:minp}) is convex w.r.t. each one of the factors $\mb{A},\mb{B}$ and $\mb{C}$ separately but not w.r.t. all of them. Moreover, due to the regularization term, it is non-separable w.r.t to each one of the matrix factors. As a result, minimizing the objective function in~(\ref{eq:minp}) alternatingly w.r.t. the BTD factors (i.e., in a BCD way with blocks the matrices $\mb{A},\mb{B}$ and $\mb{C}$) would not lead to closed-form solutions, which are desirable in an iterative algorithm. Capitalizing  on our previous work on low-rank matrix factorization \cite{grk19a}, we curb that problem by following a BSUM approach for the objective function in~(\ref{eq:minp}). The idea is that at each iteration of the BSUM scheme the BTD factors can be computed in {\it closed form} by minimizing appropriate upper bound functions of their initial objectives. Provided that these functions satisfy certain conditions~\cite{hrlp16}, the BSUM procedure is guaranteed to converge to stationary points of the objective function of the original minimization problem. 

To be more specific, using the mode-1 unfolding of $\bc{Y}$ in~(\ref{eq:minp}) and $\bc{X} = \sum_{r=1}^R\mb{A}_r\mb{B}_r^{\T}\circ\mb{c}_r$ (cf.~(\ref{eq:X1})), the objective function w.r.t. $\mb{A}$ at iteration $k$ may be expressed as follows
\begin{align}
&f_{\mb{A}}(\mb{A}|\mb{B}^k,\mb{C}^k) = \frac{1}{2}\left\|\mb{Y}_{(1)}^{\T}-\mb{P}^k\mb{A}^{\T}\right\|_{\F}^2 + \nonumber \\
&\lambda\sum_{r=1}^{R}\sqrt{\left(\sum_{l=1}^{L}\sqrt{\|\mb{a}_{rl}\|_{2}^2 + \|\mb{b}_{rl}^k\|_{2}^2 + \eta^2}\right)^2 + \|\mb{c}_r^k\|_{2}^2 + \eta^2},
\label{eq:objA}    
\end{align}
where $\mb{P}^k \triangleq \mb{B}^k\odot\mb{C}^k$. To allow this sub-problem to have closed-form solution for $\mb{A}$, we define a local tight upper bound function of~(\ref{eq:objA}) as a rough second-order Taylor approximation of $f_{\mb{A}}(\mb{A}|\mb{B}^k,\mb{C}^k)$ around $\mb{A}^k$. Namely:  
\begin{align}
& g_{\mb{A}}(\mb{A}|\mb{A}^k,\mb{B}^k,\mb{C}^k) = f_{\mb{A}}(\mb{A}^k|\mb{B}^k,\mb{C}^k) + \text{tr}\{(\mb{A}-\mb{A}^k) \nonumber \\ 
&\nabla_{\mb{A}}f_{\mb{A}}(\mb{A}^k|\mb{B}^k,\mb{C}^k)\} + \frac{1}{2}\text{vec}(\mb{A}-\mb{A}^k)^{\T}\bar{\mb{H}}_{\mb{A}^k}\text{vec}(\mb{A}-\mb{A}^k),
\label{eq:apprA}
\end{align}
where the $ILR \times ILR$ approximate Hessian matrix $\bar{\mb{H}}_{\mb{A}^k}$ of $f_{\mb{A}}(\mb{A}|\mb{B}^k,\mb{C}^k)$ at $\mb{A}^k$ is given (in analogy with~\cite{grk19a}) by 
\begin{align}
\bar{\mb{H}}_{\mb{A}^k} = \mb{I}_I\otimes(\mb{P}^{k\T}\mb{P}^k + \lambda\mb{D}^k),
\label{eq:HAk}
\end{align}
with $\mb{D}^k \triangleq (\mb{D}_1^k\otimes\mb{I}_L)\mb{D}_2^k$. $\mb{D}_1^k$ is an $R\times R$ diagonal matrix, whose $r$th diagonal entry is  
\begin{eqnarray}
\lefteqn{\mb{D}_1^k(r,r) =} \nonumber \\
& &\left[\left(\sum_{l=1}^{L}\sqrt{\|\mb{a}_{rl}^k\|_{2}^2 + \|\mb{b}_{rl}^k\|_{2}^2 +\eta^2}\right)^2 + \|\mb{c}_r^k\|_{2}^2 + \eta^2\right]^{-1/2}  
\label{eq:D1}
\end{eqnarray}
and $\mb{D}_2^k$ is an $RL\times RL$ diagonal matrix, whose $((r-1)L+l)$th diagonal entry is 
\begin{eqnarray}
\lefteqn{\mb{D}_2^k((r-1)L+l,(r-1)L+l) =} \nonumber \\
& & (\|\mb{a}_{rl}^k\|_{2}^2 + \|\mb{b}_{rl}^k\|_{2}^2 +\eta^2)^{-1/2}.
\label{eq:D2}
\end{eqnarray}
The conditions of BSUM are satisfied by the majorization function in~(\ref{eq:apprA}) if $\bar{\mb{H}}_{\mb{A}^k}$ and $\bar{\mb{H}}_{\mb{A}^k} - \mb{H}_{\mb{A}^k}$ are both positive semi-definite matrices, where $\mb{H}_{\mb{A}^k}$ is the actual Hessian of $f_{\mb{A}}(\mb{A}|\mb{B}^k,\mb{C}^k)$ at $\mb{A}^k$~\cite{hrlp16}. From~(\ref{eq:HAk}) it is obvious that $\bar{\mb{H}}_{\mb{A}^k}$ is positive semi-definite, while the positive semi-definiteness of $\bar{\mb{H}}_{\mb{A}^k} - \mb{H}_{\mb{A}^k}$ is proved in Appendix~\ref{sec:A}. Then minimizing $g_{\mb{A}}(\mb{A}|\mb{A}^k,\mb{B}^k,\mb{C}^k)$ w.r.t to $\mb{A}$ results in the following analytical expression for the estimate of $\mb{A}$ at iteration $k+1$:
\begin{align}
\mb{A}^{k+1} = \mb{Y}_{(1)}\mb{P}^k(\mb{P}^{k{\T}}\mb{P}^k+\lambda\mb{D}^{k})^{-1}. 
\label{eq:Aupdate}
\end{align}
Similarly, $\mb{B}^{k+1}$ can be found from the minimization of  $g_{\mb{B}}(\mb{B}|\mb{A}^k,\mb{B}^k,\mb{C}^k)$, which has an analogous form with~(\ref{eq:apprA}) with $\bar{\mb{H}}_{\mb{B}^k} = \mb{I}_J\otimes(\mb{Q}^{k\T}\mb{Q}^k + \lambda\mb{D}^k)$ and is a tight upper bound around $\mb{B}^k$ of 
\begin{align}
&f_{\mb{B}}(\mb{B}|\mb{A}^k,\mb{C}^k) = \frac{1}{2}\left\|\mb{Y}_{(2)}^{\T}-\mb{Q}^k\mb{B}^{\T}\right\|_{\F}^2 + \nonumber \\
&\lambda\sum_{r=1}^{R}\sqrt{\left(\sum_{l=1}^{L}\sqrt{\|\mb{a}_{rl}^k\|_{2}^2 + \|\mb{b}_{rl}\|_{2}^2 + \eta^2}\right)^2 + \|\mb{c}_r^k\|_{2}^2 + \eta^2},
\label{eq:objB}    
\end{align}
with (cf.~(\ref{eq:X2})) $\mb{Q}^k \triangleq \mb{C}^k\odot\mb{A}^k$. The unique solution of $\underset{\mb{B}}{\mathrm{min}}\; g_{\mb{B}}(\mb{B}|\mb{A}^k,\mb{B}^k,\mb{C}^k)$ is given by
\begin{align}
\mb{B}^{k+1} = \mb{Y}_{(2)}\mb{Q}^k(\mb{Q}^{k{\T}}\mb{Q}^k+\lambda\mb{D}^{k})^{-1}.    
\label{eq:Bupdate}
\end{align}
Finally, the objective function w.r.t. $\mb{C}$ may be expressed as 
\begin{align}
&f_{\mb{C}}(\mb{C}|\mb{A}^k,\mb{B}^k) = \frac{1}{2}\left\|\mb{Y}_{(3)}^{\T}-\mb{S}^k\mb{C}^{\T}\right\|_{\F}^2 + \nonumber \\
&\lambda\sum_{r=1}^{R}\sqrt{\left(\sum_{l=1}^{L}\sqrt{\|\mb{a}_{rl}^k\|_{2}^2 + \|\mb{b}_{rl}^k\|_{2}^2 + \eta^2}\right)^2 + \|\mb{c}_r\|_{2}^2 + \eta^2},
\label{eq:objC}    
\end{align}
where (cf.~(\ref{eq:X3})) 
\[
\mb{S}^k \triangleq \left[\begin{array}{ccc} (\mb{A}_1^k\odot_{\mathrm{c}} \mb{B}_1^k)\mb{1}_L & \cdots & (\mb{A}_R^k\odot_{\mathrm{c}} \mb{B}_R^k)\mb{1}_L\end{array}\right].
\]
The factor $\mb{C}^{k+1}$ is found from $\underset{\mb{C}}{\mathrm{min}}\; g_{\mb{C}}(\mb{C}|\mb{A}^k,\mb{B}^k,\mb{C}^k)$ as
\begin{align}
\mb{C}^{k+1} = \mb{Y}_{(3)}\mb{S}^k(\mb{S}^{k{\T}}\mb{S}^k+\lambda\mb{D}_1^{k})^{-1},  
\label{eq:Cupdate}
\end{align}
where the locally upper bound function $g_{\mb{C}}(\mb{C}|\mb{A}^k,\mb{B}^k,\mb{C}^k)$ has an  analogous form with $g_{\mb{A}}(\mb{A}|\mb{A}^k,\mb{B}^k,\mb{C}^k)$ in~(\ref{eq:apprA}) with 
\[
\bar{\mb{H}}_{\mb{C}^k} = \mb{I}_K\otimes(\mb{S}^{k\T}\mb{S}^k + \lambda\mb{D}_1^k).
\]

Summarizing the above, the steps of the proposed algorithm, which alternatingly solves for $\mb{A}$, $\mb{B}$, and $\mb{C}$, in that order, are tabulated as Algorithm~1. 
\begin{table}
\centering
\title{Algorithm 1: BTD-HIRLS algorithm}
 \begin{tabular}{|l|}
 \hline \\
  {\bf Algorithm 1}: BTD-HIRLS algorithm\\ \hline 
  Input: $\cal{Y}$,$\lambda,R_{\mathrm{ini}},L_{\mathrm{ini}}$ \\
  Initialize: $k=0, \mb{A}^0,\mb{B}^0, \mb{C}^0$    \\
  \bf{repeat}\\
	  \hspace{0.3cm} Compute $\mb{D}_1^k,\mb{D}_2^k$ from (\ref{eq:D1}) and (\ref{eq:D2}) \\
	   \hspace{0.3cm} $\mb{D}^k \leftarrow (\mb{D}_1^k\otimes\mb{I}_L)\mb{D}_2^k$ \\
	   \hspace{0.3cm} $\mb{P}^k \leftarrow \mb{B}^k\odot\mb{C}^k$\\
    \hspace{0.3cm} $\mb{A}^{k+1} \leftarrow \mb{Y}_{(1)}\mb{P}^k(\mb{P}^{k{\T}}\mb{P}^k+\lambda\mb{D}^{k})^{-1}$ \\
    \hspace{0.3cm} $\mb{Q}^k \leftarrow \mb{C}^k\odot\mb{A}^{k}$\\
    \hspace{0.3cm} $\mb{B}^{k+1} \leftarrow \mb{Y}_{(2)}\mb{Q}^k(\mb{Q}^{k{\T}}\mb{Q}^k+\lambda\mb{D}^{k})^{-1}$ \\
    \hspace{0.3cm} $\mb{S}^k \leftarrow \left[\begin{array}{ccc} (\mb{A}_1^{k}\odot_{\mathrm{c}} \mb{B}_1^{k})\mb{1}_L & \cdots & (\mb{A}_R^{k}\odot_{\mathrm{c}} \mb{B}_R^{k})\mb{1}_L\end{array}\right]$ \\
    \hspace{0.3cm} $\mb{C}^{k+1} \leftarrow \mb{Y}_{(3)}\mb{S}^k(\mb{S}^{k{\T}}\mb{S}^k+\lambda\mb{D}_1^{k})^{-1}$ \\
    \hspace{0.4cm}$k \leftarrow k+1$\\
    \bf{until} {\it convergence} \\
    \hline
 \end{tabular}
\end{table}
As explained in Appendix B, the proposed algorithm is a sort of  {\it hierarchical} iterative reweighted least squares (HIRLS) scheme, fully adjusted to promote block and column sparsity in the BTD model. This may be also seen from the expressions of $\mb{A}^{k+1}, \mb{B}^{k+1}$, and $\mb{C}^{k+1}$ given above and the form of the diagonal weighting matrices $\mb{D}_1^{k}$ and $\mb{D}_2^k$.  
Indeed, if $R$ and $L$ are overestimated, reweighting via $\mb{D}_1$ imposes \emph{jointly} block sparsity on $\mb{A}$ and $\mb{B}$ and column sparsity on $\mb{C}$, hence helping in estimating $R$. In addition, reweighting via $\mb{D}_2$ promotes column sparsity \emph{jointly} on the corresponding blocks of $\mb{A}$ and $\mb{B}$, thus estimating $L_r$'s. This mechanism, combined with an appropriate selection of $\lambda$, can reveal the actual value of $R$ and the true block-term ranks $L_r$'s, as it is also empirically demonstrated in the next section.

It should be noted that the majorization functions $g_{\mb{A}}$, $g_{\mb{B}}$, and $g_{\mb{C}}$ used previously for the derivation of the proposed algorithm are {\it quadratic upper bound} functions that satisfy Assumption~A~\cite[Table~3]{hrlp16} required for BSUM. In addition, minimization of these functions w.r.t. the BTD factors at each iteration of the algorithm leads in all cases to unique solutions. Hence, according to Theorem~1 of~\cite{hrlp16}, every limit point of the BTD-HIRLS algorithm is a stationary point of the initial objective function~(\ref{eq:minp}).

A notable feature of the proposed algorithm is that the closed-form expressions for the BTD factors comprise matrix operations only and relatively small-size matrix inversions, which is translated to relatively low computational complexity. In contrast, in~\cite{gofzc20}, the BTD factors are not computed in closed form but via a group-sparsity promoting iterative procedure in each iteration of the algorithm, which results in an increase of the computation time. Further reduction in the computational complexity of the proposed algorithm is possible by eliminating negligible columns (column pruning) in the course of the iterations (as in~\cite{grk19a}). 

\section{Simulation Results}
\label{sec:examples}

In this section, we report indicative simulation results with both synthetic and real data for evaluating the performance of the proposed algorithm. For comparison purposes, the classical BTD-ALS algorithm of~\cite{ldl08c}, which makes no use of any low-rank regularization, and the BTD-AGL algorithm of~\cite{gofzc20}, which minimizes the objective function defined in~(\ref{eq:ABC}) enhanced by a proximal term, are also tested. It should be noted that a block-pruning mechanism is implemented in both BTD-HIRLS and BTD-AGL, that is, $\mb{A}_r,\mb{B}_r$ blocks that correspond to columns of $\mb{C}$ with negligible energy are removed as the algorithms progress. This can be applied in both of the aforementioned algorithms due to their group-sparsity imposing characteristics. 

\subsection{Synthetic data}

In all cases, we generate BTD tensors $\bc{X}$ contaminated by additive noise, i.e., 
$\bc{Y}=\bc{X}+\sigma\bc{N}$, 
where $\bc{N}$ contains zero-mean, independent and identically distributed (i.i.d) Gaussian entries of unit variance and $\sigma$ is set so that we get a given signal-to-noise ratio (SNR), with SNR in~dB defined as $\mathrm{SNR}=10\log_{10}\|\bc{X}\|_{\F}^2/(\sigma^2\|\bc{N}\|_{\F}^2)$. The entries of the matrices $\mb{A}_r$ and $\mb{B}_r$ and the vectors $\mb{c}_r$ have been also sampled from i.i.d. zero-mean Gaussian distributions of unit variance. The tensor approximation is measured with the normalized mean squared error (NMSE) over the blocks, defined as $\mathrm{NMSE}(\hat{\mathbf{A}},\hat{\mathbf{B}},\hat{\mathbf{C}})= \frac{1}{R}\sum^R_{r=1} \frac{\|\mathbf{A}_r\mathbf{B}^{\T}_r\circ \mathbf{c}_r - \hat{\mathbf{A}}_r\hat{\mathbf{B}}^{\T}_r\circ \hat{\mathbf{c}}_r\|_{\F}^2}{\|\mathbf{A}_r\mathbf{B}^{\T}_r\circ \mathbf{c}_r\|_{\F}^2}$, where $(\mathbf{A,B,C})$ and $(\hat{\mathbf{A}},\hat{\mathbf{B}},\hat{\mathbf{C}})$ denote the true and the estimated tensor factors, respectively. To calculate this metric, a linear assignment problem is solved to resolve the permutation ambiguity.\footnote{The Matlab 2019b {\tt matchpairs} function was employed for this purpose.} When $R$ is overestimated (as in BTD-HIRLS and BTD-AGL), the NMSE over blocks is calculated on the basis of those of the $R$ block terms that are ``closer" to the true ones. For the stopping criterion we use the relative difference between two consecutive values of the reconstruction error. The algorithms stop either when the relative difference becomes less than $10^{-5}$ or a maximum of 200~iterations is reached. The regularization parameter $\lambda$ was empirically observed to depend on the dimensions of the tensor and the model ranks as well as on the noise strength. Hence, for the selection of the value of $\lambda$ in BTD-HIRLS we employ the heuristic rule $\lambda=L_{\mathrm{ini}}R_{\mathrm{ini}}(I+J+K)\hat{\sigma}$ with $\hat{\sigma}$ being a guess of the standard deviation of the noise. 
Initializations are performed randomly.
The $\gamma$-sweeping procedure employed in~\cite{gofzc20} is also adopted here for BTD-AGL, namely, for each single case (corresponding to a given realization and an initialization) BTD-AGL is applied five times with five different (and increasing) $\gamma$'s (as in~\cite{gofzc20}) and the estimates from each run are used to initialize the next one.

\subsubsection{Performance in the presence of noise}

First, we test the algorithms at different SNR values. We set $I=60$, $J=50$ and $K=55$. The true $R$ is set to~5 and the $L_r$'s are integer numbers chosen uniformly at random from the set $\{2,3,4,5,6,7,8,9\}$. The noisy tensors are generated as described above. Since the ranks of the model are in general unknown, we initialized BTD-HIRLS and BTD-AGL with overestimates of the true ones, namely ${R}_{\mathrm{ini}}=10$ and $L_{\mathrm{ini}}=10$ for all factor blocks. For BTD-ALS it was assumed that the true $R$ is known,  while all $L_r$'s were overestimated to 10. All algorithms were randomly initialized 10~times and their best run, in terms of the NMSE, was kept. In Table~\ref{tab:noise_exp} we report the median NMSEs of the results obtained over 10~independent realizations of the experiment. The average run-times (with Matlab~2019b in a MacBook~Pro, 2.6~GHz 6-Core Intel Core i7, 16~GB 2667~MHz DDR4) are also reported. 
\begin{table}
    \centering
        \caption{NMSE and run-time comparison of BTD-HIRLS, BTD-AGL and BTD-ALS for different SNR values.}
    \begin{tabular}{c|c|c|c|c|c|} 
    \cline{2-6} 
   & \multicolumn{4}{|c|}{SNR (dB)}  & \multirow{2}{*}{aver. time (sec)} \\ \cline{2-5}
     & \multicolumn{1}{c|}{5} & \multicolumn{1}{c|}{10} & \multicolumn{1}{c|}{15}  & \multicolumn{1}{c|}{20} & \\ \hline
      \multicolumn{1}{|c|}{BTD-HIRLS }  & 0.0128   & 0.0040 &  0.0012 & 0.0004 & 2.5\\ \hline
   \multicolumn{1}{|c|}{BTD-AGL}  & 0.0197  & 0.0053 &  0.0015 & 0.0004 & 68 \\ \hline 
   \multicolumn{1}{|c|}{BTD-ALS}  & 0.0893  & 0.0562 & 0.0561 & 0.0561 & 2.2\\ \hline
    \end{tabular}
    \label{tab:noise_exp}
\end{table}
The proposed method is seen to outperform both BTD-ALS and BTD-AGL in terms of NMSE. Moreover, BTD-AGL is considerably more computationally costly. Note that this should not be attributed to the $\gamma$-sweeping procedure only. Even a single run of BTD-AGL takes more time because, in contrast to BTD-HIRLS (and BTD-ALS), BTD-AGL does not rely on closed-form solutions for the updates of $\mb{A}$, $\mb{B}$, and $\mb{C}$ but it instead involves  separate iterative procedures for the solution of each of the sub-problems. 

\begin{figure}
\centerline{\includegraphics[width=0.45\textwidth,height=0.3\textwidth]{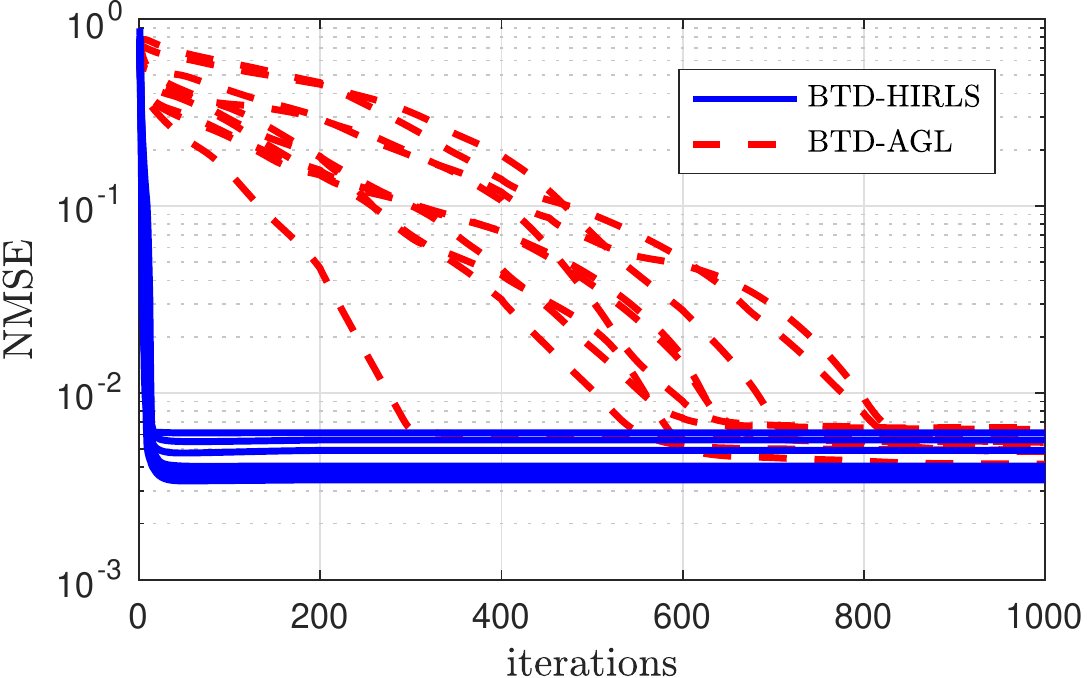}}
\caption{NMSE of BTD-HIRLS and BTD-AGL vs. iterations for SNR=10~dB.}
\label{fig:nmse_vs_iters}
\end{figure}

Furthermore, BTD-HIRLS exhibits a much higher rate of convergence than BTD-AGL, as demonstrated in Fig.~\ref{fig:nmse_vs_iters}, where the evolution of NMSE's is plotted versus the number of iterations, for SNR=10~dB and all 10~realizations of the experiment.
To facilitate the comparison of the 200 iterations of BTD-HIRLS with the BTD-AGL sweeping runs, its curves have been extended all the way to 1000~iterations based on the NMSE value at iteration~200. It should be clear from Fig.~\ref{fig:nmse_vs_iters} that the $\gamma$-sweeping procedure is indeed necessary for BTD-AGL as in most cases it does not converge before the 3rd~round of it, that is, 600~iterations. In contrast, the proposed method converges fast, requiring no more than 50~iterations in all realizations.

\subsubsection{Success rates for the recovery of $R$ and $L_r$'s}

\begin{figure*}
\begin{tabular}{c c}
{\includegraphics[width=0.45\textwidth,height=0.4\textwidth]{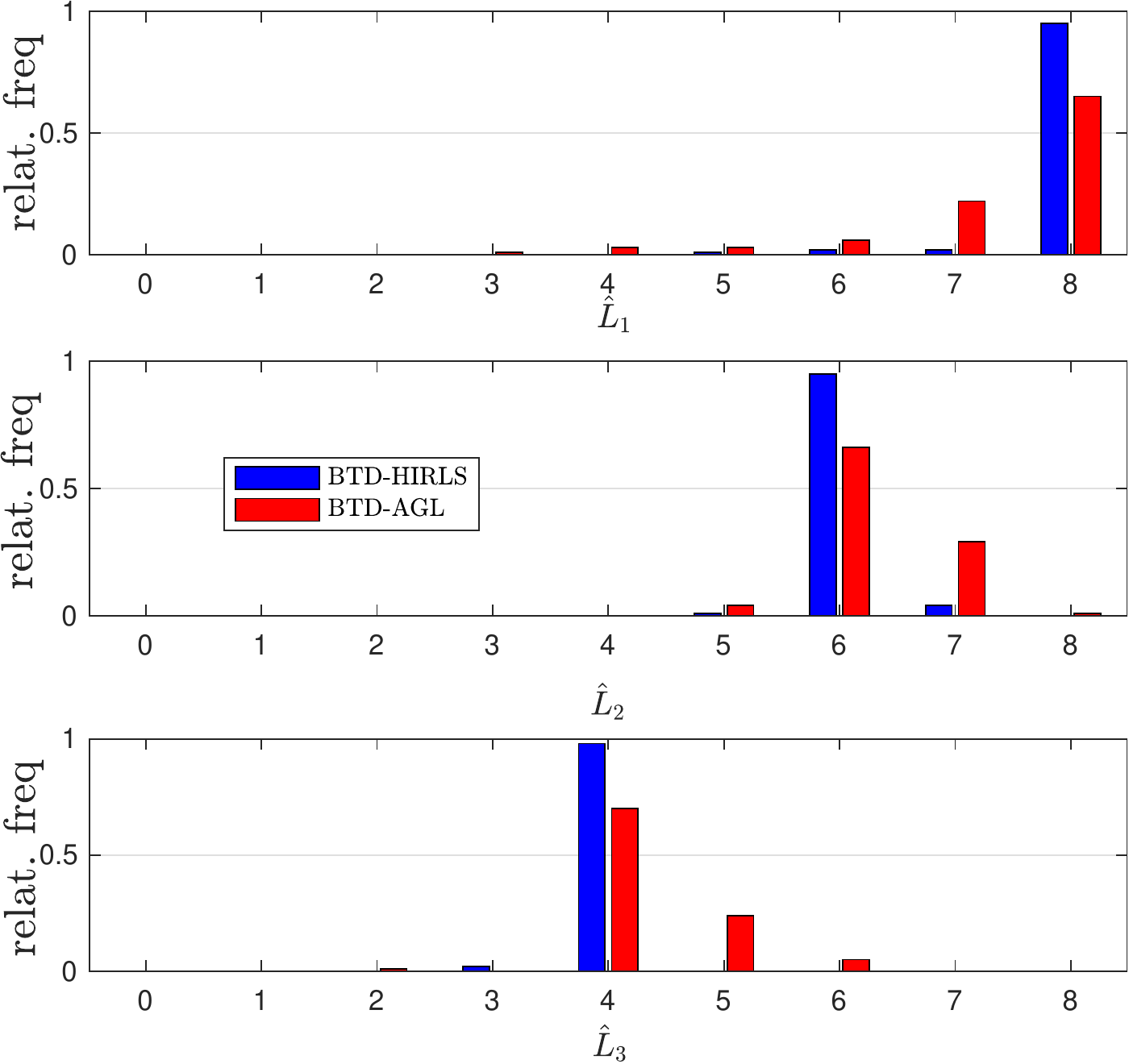}} & {\includegraphics[width=0.45\textwidth,height=0.4\textwidth]{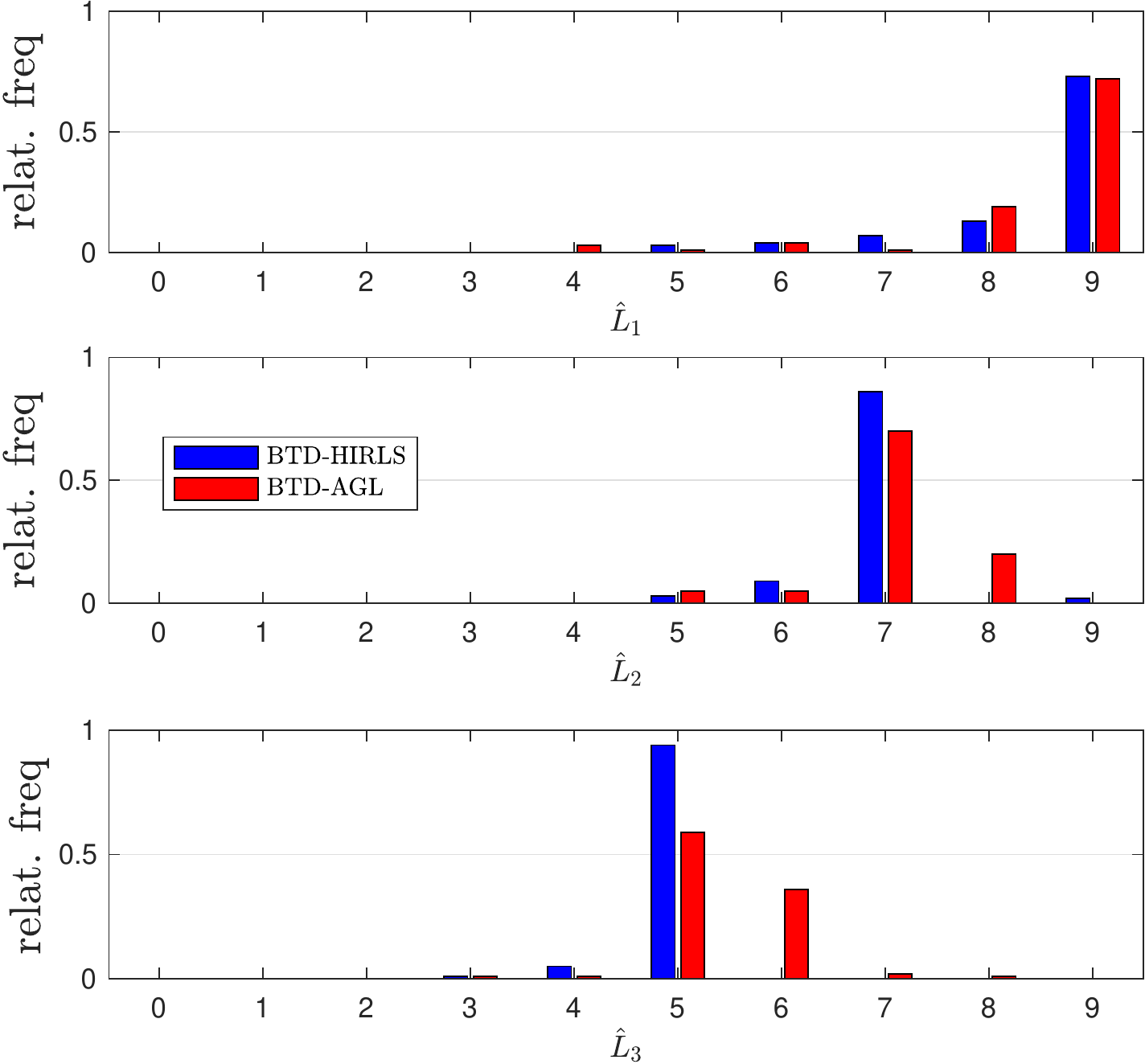}}  \\
(a) & (b) \\
{\includegraphics[width=0.35\textwidth,height=0.15\textwidth]{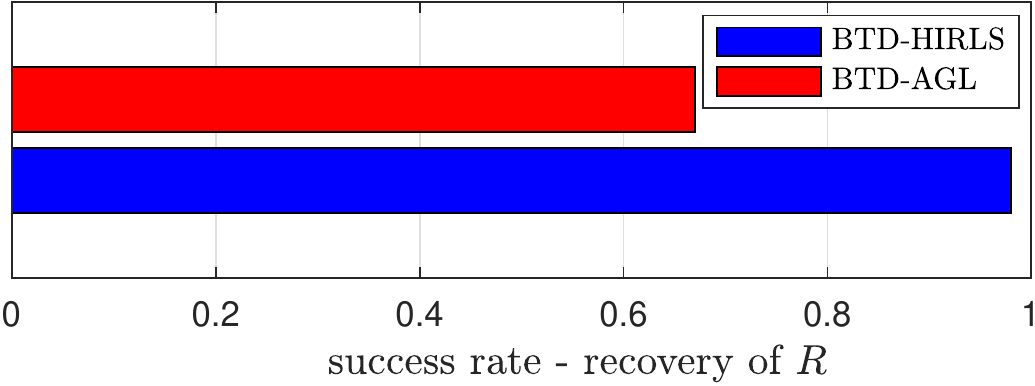}}& {\includegraphics[width=0.35\textwidth,height=0.15\textwidth]{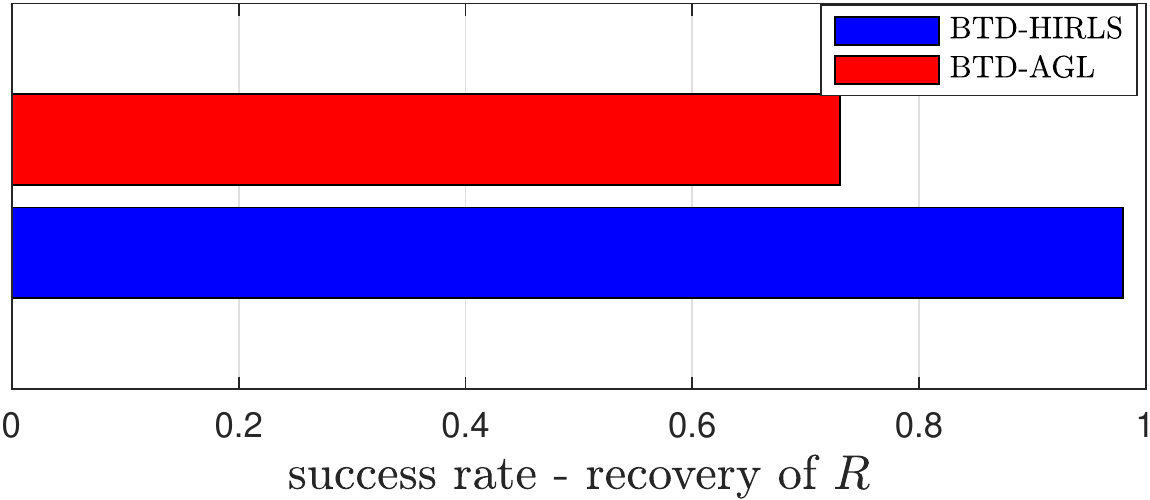}}\\
(c) & (d)
\end{tabular}
\caption{Relative frequencies of the estimated values of $L_r$'s and success rate (\%) of estimating $R$ via BTD-HIRLS and BTD-AGL. SNR=15~dB. a) Scenario~I: $\sum^R_{r=1}L_r \leq \mathrm{min}(I,J)$; true $L_r$'s are $L_1=8, L_2=6, L_3=4$ b) Scenario~II: $\sum^R_{r=1}L_r > \mathrm{min}(I,J)$; true $L_r$'s are $L_1=9, L_2=7, L_3=5$ (c) success rate of estimating $R$ for Scenario~I and (d) same for Scenario~II.}
\label{fig:succ_rate}
\end{figure*}

In this part, our aim is to demonstrate the capability of BTD-HIRLS to revealing the true model structure. We set SNR=15~dB and we estimate the success rates in the estimation of $R$ as well as of the $L_r$'s for 100~different realizations of the experiment. Here tensors are of dimensions $18\times 18\times 10$ and the true number of the blocks, $R$, is set to~3. Again, we compare BTD-HIRLS with BTD-AGL and for both algorithms we over-estimate $R$ and $L_r$'s as $R_{\mathrm{ini}}=10$ and $L_{\mathrm{ini}} = 10$ for all block terms. We examine two different scenarios:\\
a) {\it Scenario I}: The true block ranks are $L_1=8$, $L_2=6$ and $L_3=4$. In this case, $\sum^R_{r=1}L_r =\mathrm{min}(I,J)$, that is, the sufficient uniqueness condition is met. As it can be seen in Fig.~\ref{fig:succ_rate}(a), BTD-HIRLS achieves success rates higher than 90\% for all $L_r$'s outperforming BTD-AGL in the task of revealing the true $L_r$'s. This can be explained by the properties of the regularizer of BTD-HIRLS which is carefully designed so as to better capture the structure of the decomposition model. The latter is also verified in  Fig.~\ref{fig:succ_rate}(c), where we can see that BTD-HIRLS is more efficient than BTD-AGL when it comes to the success in estimating the number of block terms, $R$. \\
b) {\it Scenario II}: In this more challenging setting, we set $L_1=9$, $L_2=7$, and $L_3=5$. Thus, we now have $\sum^R_{r=1}L_r > \mathrm{min}(I,J)$ and hence the sufficient uniqueness condition is violated. However, it can be observed in Fig.~\ref{fig:succ_rate}(b) that BTD-HIRLS reveals all $L_r$'s with high relative frequencies, slightly outperforming BTD-AGL. Moreover, the success rate of accurately estimating $R$ remains high as in the case of the ``benign" scenario (cf. Fig.~\ref{fig:succ_rate}(d)). 

\subsection{Experimenting with real data: Hyper-spectral image de-noising} 

Hyper-spectral images are known to exhibit high coherence both in the spectral and the spatial domain~\cite{grk19a}. As a result, low-rank matrix and tensor factorization methods have been widely employed to address related problems such as restoration, super resolution, de-noising, etc.~\cite{grk19a,zfhw19,xzq19}. Here, we consider the recovery of a hyper-spectral image from its noise-corrupted version with additive Gaussian noise of SNR=5~dB, with the aid of the BTD-HIRLS and BTD-AGL algorithms. The idea is to recover the image as a low-rank BTD approximation of the noisy one, capitalizing on the low-rank structure of the HSI, which allows the removal of the high-rank noise~\cite{grk19a}. In both cases, we set $R_{\mathrm{ini}}=50$ and $L_{\mathrm{ini}}=10$. Note that, for HSI data, and using the HSI spectral unmixing jargon, $R$ is related to the number of the end-members (i.e., spectral signatures of the materials that exist in the depicted scene) while the $L_r$'s reflect the ranks of the corresponding $R$ abundance maps (i.e., the images of the percentages of a given material in the given image). 
In this example, we consider the Washington DC Mall AVIRIS image captured at $m = 191$ contiguous spectral bands in the~0.4 to~2.4 $\mu m$ region of the visible and infrared spectrum~\cite{grk19a}. The size of the image is $150\times 150$ pixels. 

\begin{figure}
{\includegraphics[width=0.45\textwidth,height=0.3\textwidth]{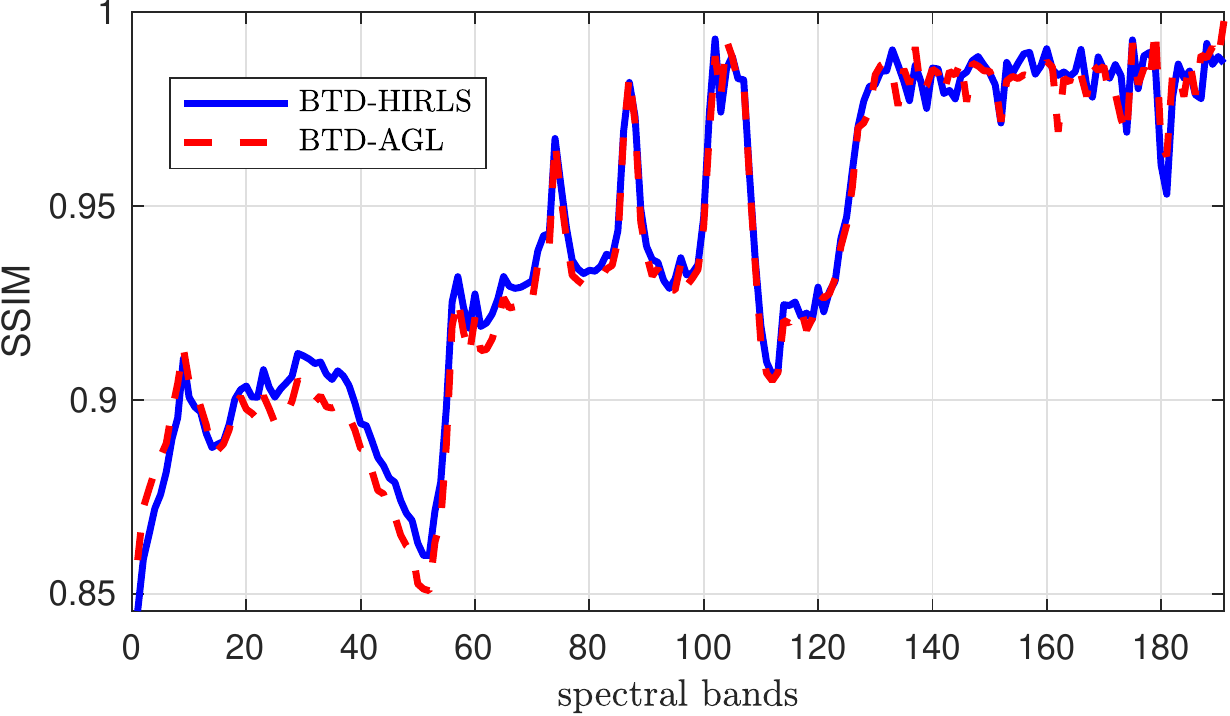}}\\
\caption{SSIM of the hyper-spectral images recovered via BTD-HIRLS and BTD-AGL.}
\label{fig:ssim}
\end{figure}

\begin{figure}
    \begin{tabular}{c c}
        \includegraphics[width=0.22\textwidth,height=0.2\textwidth]{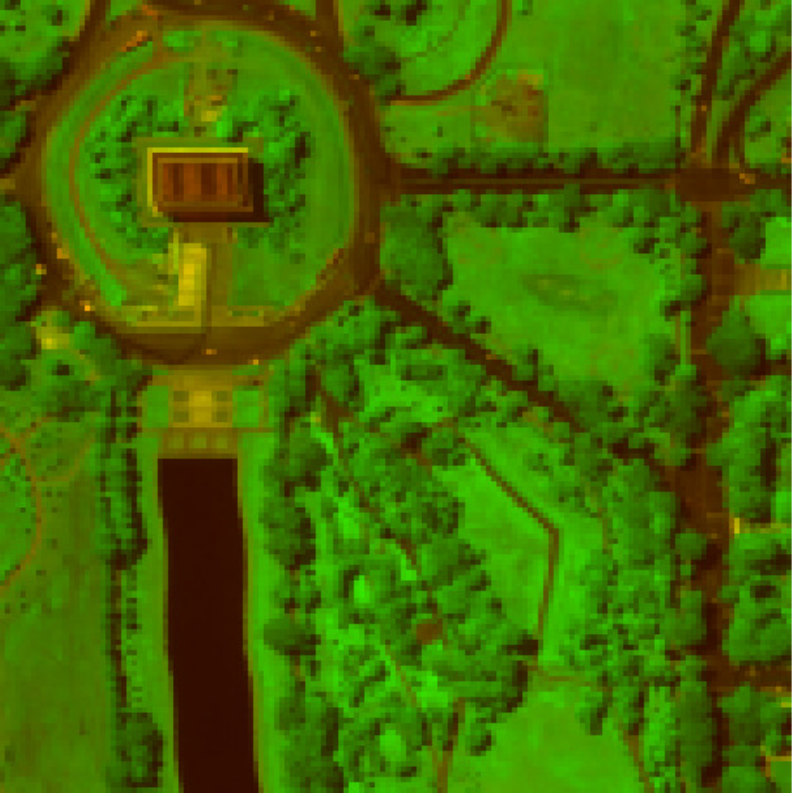} &  \includegraphics[width=0.22\textwidth,height=0.2\textwidth]{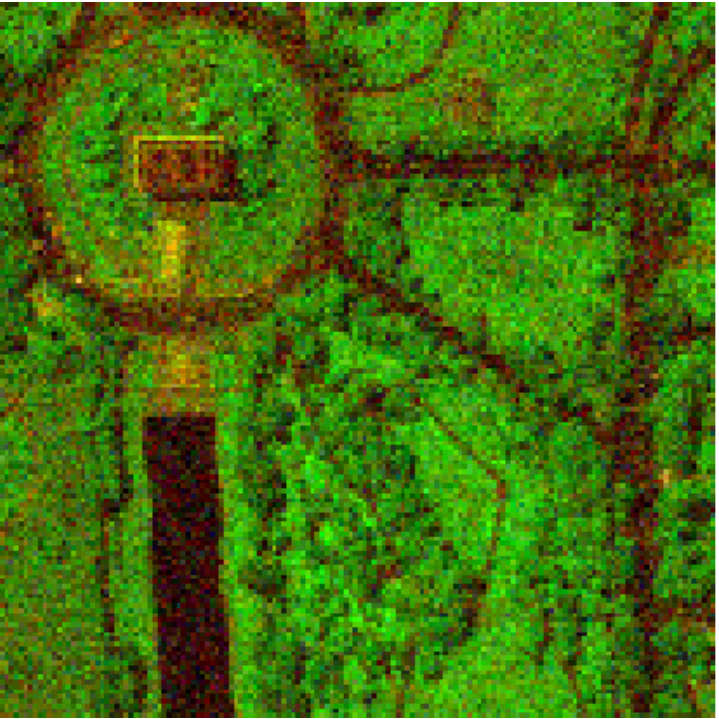}   \\
        (a) & (b) \\
        \includegraphics[width=0.22\textwidth,height=0.2\textwidth]{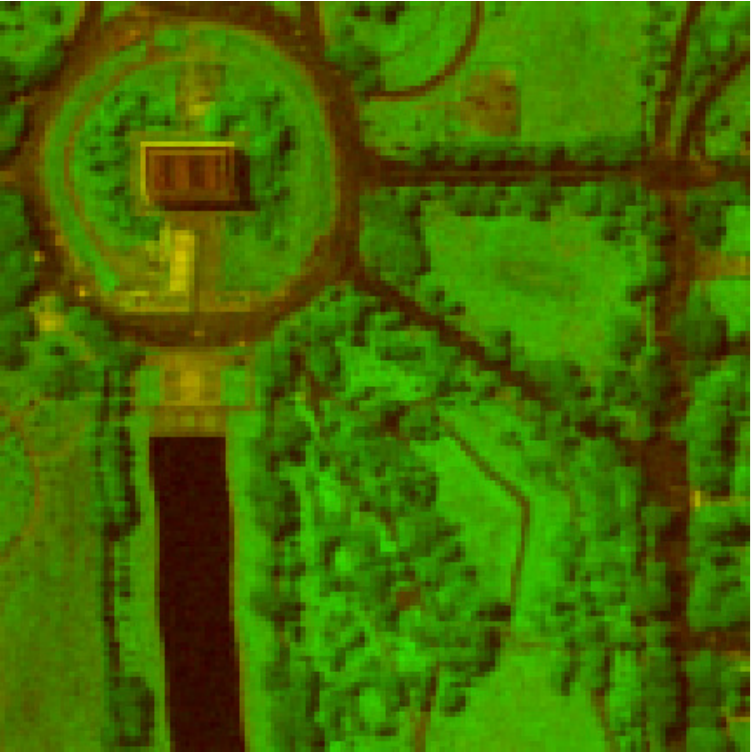} & \includegraphics[width=0.22\textwidth,height=0.2\textwidth]{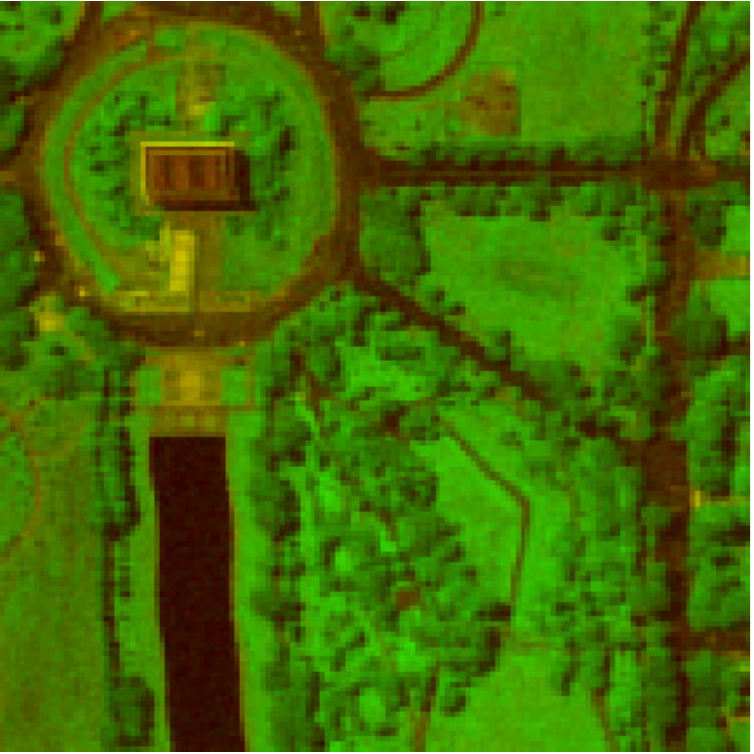} \\
        (c) & (d) 
    \end{tabular}
    \caption{False RGB color images of Washington DC Mall AVIRIS hyper-spectral image. (a) Original (b) Noisy (c) De-noised with BTD-AGL (d) De-noised with BTD-HIRLS.}
    \label{fig:false_rgb_hyperspectal}
\end{figure}

Fig.~\ref{fig:ssim} plots the values of the structural similarity index (SSIM) between the original and the de-noised images, while the results can be visually inspected in Fig.~\ref{fig:false_rgb_hyperspectal}, where RGB false color images reconstructed from bands (24,64,135) are shown. 
The two algorithms are seen to perform equally well in this experiment. Both estimate $R$ as~8, which agrees with the true number of the end-members in the scene depicted~\cite{grk19a}. 

\section{Conclusions}
\label{sec:concls}

The challenging problem of efficiently and effectively estimating the model structure and parameters of a BTD has recently received special attention due to the increasing application range of this tensor model. This paper briefly reviews the related literature and reports our recent results on this topic, which are based on an appropriate extension to the BTD model of our earlier rank-revealing work on low-rank matrix and tensor approximation. The idea is to impose column sparsity \emph{jointly} on the factors and in a \emph{hierarchical} manner that matches the structure of the model, and successively estimate the ranks as the numbers of factor columns of non-negligible magnitude, with the aid of alternating hierarchical IRLS. The proposed method enjoys fast convergence and low computational complexity, also allowing the negligible columns to be pruned in the course of the procedure. Simulation results that demonstrate the effectiveness of our method in accurately estimating both the ranks and the factors in both synthetic and real-world scenarios are reported.

Future work will include the development of constrained variants of the method and (semi-)automatic ways of tuning its regularization parameter. 

\appendices

\section{Proof of Positive Semi-definiteness}
\label{sec:A}

Let the $ILR \times 1$ vector $\mb{a} \triangleq \mathrm{vec}(\mb{A})$ be $[a_{11},\ldots,a_{1d},\ldots,a_{I1},\ldots,a_{Id}]^{\T}$ where $d\triangleq LR$. 
After some tedious algebra, it can be shown that $\bar{\mb{H}}_{\mb{A}}-\mb{H}_{\mb{A}} = \lambda \mb{U}$, where $\mb{U}$ is a $Id \times Id$ matrix which consists of $I^2$ $d \times d$ diagonal blocks denoted as $\mb{U}_{ij}$ with $i,j = 1,2,\ldots,I$ and shown in~(\ref{eq:Ukj}) (see top of the next page).
\begin{figure*}
\begin{eqnarray}
\mb{U}_{ij} = \mathrm{diag}\Big(&\frac{a_{i1}a_{j1}\left[\left(\|\mb{a}_{11}\|_2^2+\|\mb{b}_{11}\|_2^2\right)^{1/2} + \left(\sum_{l=1}^L\sqrt{\|\mb{a}_{1l}\|_2^2+\|\mb{b}_{1l}\|_2^2}\right)^2 + \|\mb{c}_1\|_2^2 \right]}{\left[\left(\sum_{l=1}^L\sqrt{\|\mb{a}_{1l}\|_2^2+\|\mb{b}_{1l}\|_2^2}\right)^2 + \|\mb{c}_1\|_2^2\right]^{3/2}\left(\|\mb{a}_{11}\|_2^2+\|\mb{b}_{11}\|_2^2 \right)^{3/2}}, \ldots \nonumber \\
& \frac{a_{iL}a_{jL}\left[\left(\|\mb{a}_{1L}\|_2^2+\|\mb{b}_{1L}\|_2^2\right)^{1/2} + \left(\sum_{l=1}^L\sqrt{\|\mb{a}_{1l}\|_2^2+\|\mb{b}_{1l}\|_2^2}\right)^2 + \|\mb{c}_1\|_2^2 \right]}{\left[\left(\sum_{l=1}^L\sqrt{\|\mb{a}_{1l}\|_2^2+\|\mb{b}_{1l}\|_2^2}\right)^2 + \|\mb{c}_1\|_2^2\right]^{3/2}\left(\|\mb{a}_{1L}\|_2^2+\|\mb{b}_{1L}\|_2^2  \right)^{3/2}}, \nonumber \\ 
& \frac{a_{i(L+1)}a_{j(L+1)}\left[\left(\|\mb{a}_{21}\|_2^2+\|\mb{b}_{21}\|_2^2\right)^{1/2} + \left(\sum_{l=1}^L\sqrt{\|\mb{a}_{2l}\|_2^2+\|\mb{b}_{2l}\|_2^2}\right)^2 + \|\mb{c}_2\|_2^2 \right]}{\left[\left(\sum_{l=1}^L\sqrt{\|\mb{a}_{2l}\|_2^2+\|\mb{b}_{2l}\|_2^2}\right)^2 + \|\mb{c}_2\|_2^2\right]^{3/2}\left(\|\mb{a}_{21}\|_2^2+\|\mb{b}_{21}\|_2^2  \right)^{3/2}}, \ldots \nonumber \\ 
& \frac{a_{id}a_{jd}\left[\left(\|\mb{a}_{RL}\|_2^2+\|\mb{b}_{RL}\|_2^2\right)^{1/2} + \left(\sum_{l=1}^L\sqrt{\|\mb{a}_{Rl}\|_2^2+\|\mb{b}_{Rl}\|_2^2}\right)^2 + \|\mb{c}_R\|_2^2 \right]}{\left[\left(\sum_{l=1}^L\sqrt{\|\mb{a}_{Rl}\|_2^2+\|\mb{b}_{Rl}\|_2^2}\right)^2 + \|\mb{c}_R\|_2^2\right]^{3/2}\left(\|\mb{a}_{RL}\|_2^2+\|\mb{b}_{RL}\|_2^2  \right)^{3/2}} \Big)
\label{eq:Ukj}
\end{eqnarray}
\hrule
\end{figure*}
From~(\ref{eq:Ukj}) it follows that $\mb{U}$ can be written in the form $\mb{U} = \tilde{\mb{U}}^{\T}\tilde{\mb{U}}$, hence it is positive semi-definite. In fact, $\tilde{\mb{U}}$ is a $d \times Id$ matrix that comprises $I$ $d \times d$ diagonal blocks $\tilde{\mb{U}}_i$ which can be written as in~(\ref{eq:Uk}) (see top of the next page). 
\begin{figure*}
\begin{eqnarray}
\tilde{\mb{U}}_{i} = \mathrm{diag}\Big(&\frac{a_{i1}\left[\left(\|\mb{a}_{11}\|_2^2+\|\mb{b}_{11}\|_2^2\right)^{1/2} + \left(\sum_{l=1}^L\sqrt{\|\mb{a}_{1l}\|_2^2+\|\mb{b}_{1l}\|_2^2}\right)^2 + \|\mb{c}_1\|_2^2 \right]^{1/2}}{\left[\left(\sum_{l=1}^L\sqrt{\|\mb{a}_{1l}\|_2^2+\|\mb{b}_{1l}\|_2^2}\right)^2 + \|\mb{c}_1\|_2^2\right]^{3/4}\left(\|\mb{a}_{11}\|_2^2+\|\mb{b}_{11}\|_2^2 \right)^{3/4}}, \ldots \nonumber \\
& \frac{a_{iL}\left[\left(\|\mb{a}_{1L}\|_2^2+\|\mb{b}_{1L}\|_2^2\right)^{1/2} + \left(\sum_{l=1}^L\sqrt{\|\mb{a}_{1l}\|_2^2+\|\mb{b}_{1l}\|_2^2}\right)^2 + \|\mb{c}_1\|_2^2 \right]^{1/2}}{\left[\left(\sum_{l=1}^L\sqrt{\|\mb{a}_{1l}\|_2^2+\|\mb{b}_{1l}\|_2^2}\right)^2 + \|\mb{c}_1\|_2^2\right]^{3/4}\left(\|\mb{a}_{1L}\|_2^2+\|\mb{b}_{1L}\|_2^2  \right)^{3/4}}, \nonumber \\ 
& \frac{a_{i(L+1)}\left[\left(\|\mb{a}_{21}\|_2^2+\|\mb{b}_{21}\|_2^2\right)^{1/2} + \left(\sum_{l=1}^L\sqrt{\|\mb{a}_{2l}\|_2^2+\|\mb{b}_{2l}\|_2^2}\right)^2 + \|\mb{c}_2\|_2^2 \right]^{1/2}}{\left[\left(\sum_{l=1}^L\sqrt{\|\mb{a}_{2l}\|_2^2+\|\mb{b}_{2l}\|_2^2}\right)^2 + \|\mb{c}_2\|_2^2\right]^{3/4}\left(\|\mb{a}_{21}\|_2^2+\|\mb{b}_{21}\|_2^2  \right)^{3/4}}, \ldots \nonumber \\ 
& \frac{a_{id}\left[\left(\|\mb{a}_{RL}\|_2^2+\|\mb{b}_{RL}\|_2^2\right)^{1/2} + \left(\sum_{l=1}^L\sqrt{\|\mb{a}_{Rl}\|_2^2+\|\mb{b}_{Rl}\|_2^2}\right)^2 + \|\mb{c}_R\|_2^2 \right]^{1/2}}{\left[\left(\sum_{l=1}^L\sqrt{\|\mb{a}_{Rl}\|_2^2+\|\mb{b}_{Rl}\|_2^2}\right)^2 + \|\mb{c}_R\|_2^2\right]^{3/4}\left(\|\mb{a}_{RL}\|_2^2+\|\mb{b}_{RL}\|_2^2  \right)^{3/4}} \Big)
\label{eq:Uk}
\end{eqnarray}
\hrule
\end{figure*}
It thus follows that, since $\lambda > 0$, $\bar{\mb{H}}_{\mb{A}}-\mb{H}_{\mb{A}}$ is also positive semi-definite, which completes the proof. Analogous results for the $\mb{B}$ and $\mb{C}$ sub-problems can be similarly arrived at.

\section{The Hierarchical IRLS Nature of BTD-HIRLS}
\label{sec:B}

If $\mb{x} = \left[\begin{array}{cccc} x_1 & x_2 & \cdots & x_n\end{array}\right]^{\T}$ is a sparse vector, the IRLS algorithm for estimating $\mb{x}$ subject to an $\ell_2$ proximity criterion is derived by solving the following minimization problem at iteration $k+1$~\cite{ddfg09}:
\begin{align}
\mb{x}^{k+1} =  \arg\underset{\mb{x}}{\mathrm{min}}\frac{1}{2} \left\|\mb{b}-\mb{R}\mb{x}\right\|_{2}^2  + \frac{\lambda}{2}\sum_{i=1}^n\frac{x_i^2}{\sqrt{x_i^{k2}+\eta^2}}.
\label{eq:minx}
\end{align}
This problem admits the closed-form solution $\mb{x}^{k+1} = (\mb{R}^{\T}\mb{R} + \lambda \mb{W}^k)^{-1}\mb{R}^{\T}\mb{b}$, where $\mb{W}^k$ is a diagonal weighting matrix whose $i$th diagonal entry is given by $\mb{W}^k(i,i) = (x_i^{k2}+\eta^2)^{-1/2}$.

In the same vein, it can be shown that the closed-form expression for the BTD factor $\mb{A}^{k+1}$ given in~(\ref{eq:Aupdate}) can be also obtained by solving the following minimization problem
\begin{align}
&\mb{A}^{k+1} = \arg\underset{\mb{A}}{\mathrm{min}}\frac{1}{2}\left\|\mb{Y}_{(1)}^{\T}-\mb{P}^k\mb{A}^{\T}\right\|_{\F}^2 + \nonumber \\
 &\frac{\lambda}{2}\sum_{r=1}^{R}\frac{\left(\frac{1}{2}\sum_{l=1}^{L}\frac{\|\mb{a}_{rl}\|_{2}^2 + \|\mb{b}_{rl}^k\|_{2}^2 + \eta^2}{\sqrt{\|\mb{a}_{rl}^k\|_{2}^2 + \|\mb{b}_{rl}^k\|_{2}^2 + \eta^2}}\right)^2 + \|\mb{c}_r^k\|_{2}^2 +\eta^2}{\sqrt{\left(\sum_{l=1}^{L}\sqrt{\|\mb{a}_{rl}^k\|_{2}^2 + \|\mb{b}_{rl}^k\|_{2}^2 +\eta^2}\right)^2 + \|\mb{c}_r^k\|_{2}^2 + \eta^2}}.
\label{eq:Ak}   
\end{align}
A similar minimization problem can be defined for computing the factor $\mb{B}^{k+1}$ as in~(\ref{eq:Bupdate}), while we can also get~(\ref{eq:Cupdate}) from 
\begin{align}
&\mb{C}^{k+1} = \arg\underset{\mb{C}}{\mathrm{min}}\frac{1}{2}\left\|\mb{Y}_{(3)}^{\T}-\mb{S}^k\mb{C}^{\T}\right\|_{\F}^2 + \nonumber \\
 &\frac{\lambda}{2}\sum_{r=1}^{R}\frac{\left(\sum_{l=1}^{L}\sqrt{\|\mb{a}_{rl}^k\|_{2}^2 + \|\mb{b}_{rl}^k\|_{2}^2 + \eta^2}\right)^2 + \|\mb{c}_r\|_{2}^2 +\eta^2}{\sqrt{\left(\sum_{l=1}^{L}\sqrt{\|\mb{a}_{rl}^k\|_{2}^2 + \|\mb{b}_{rl}^k\|_{2}^2 +\eta^2}\right)^2 + \|\mb{c}_r^k\|_{2}^2 + \eta^2}}
\label{eq:Ck} 
\end{align}
By carefully inspecting the objective functions in~(\ref{eq:Ak}) and~(\ref{eq:Ck}) and comparing with the conventional IRLS objective function in~(\ref{eq:minx}), we easily recognize a hierarchical IRLS structure consisting of two separate reweighting least squares steps. Each step gives rise to a separate reweighting matrix. Namely, the first one ($\mb{D}_1$) is composed of the inverses of the outer summation terms of the regularizer in~(\ref{eq:minp}) and jointly weighs the blocks of $\mb{A},\mb{B}$, i.e. the $\mathbf{A}_r$'s and $\mathbf{B}_r$'s, and the respective columns of $\mb{C}$. The second reweighting matrix ($\mb{D}_2$) contains the inverses of the terms of the inner summation in~(\ref{eq:minp}) and jointly balances the corresponding columns of the $\mb{A}_r$'s and $\mb{B}_r$'s. It thus follows that the proposed two-level $\ell_{1,2}$ regularization naturally leads to a IRLS scheme with a corresponding hierarchy. 

\bibliographystyle{IEEEtran}
\bibliography{IEEEabrv,refs}

\begin{IEEEbiography}{Athanasios A. Rontogiannis} received the (5-yr) Diploma degree in electrical engineering from the National Technical University of Athens, Athens, Greece, in 1991, the M.A.Sc. degree in electrical and computer engineering from the University of Victoria, Victoria, BC, Canada, in 1993, and the Ph.D. degree in signal processing and communications from the National and Kapodistrian University of Athens, Athens, Greece, in 1997. From 1998 to 2003, he was a Lecturer with the University of Ioannina, Ioannina, Greece. In 2003 he joined the National Observatory of Athens, where he is currently a Research Director with the Institute for Astronomy, Astrophysics, Space Applications and Remote Sensing (IAASARS). His research interests include the general area of statistical signal and image processing with emphasis on adaptive signal processing, hyper-spectral image processing, compressive sensing, sparse and low-rank signal representations, and fast signal processing algorithms. On these topics, he has co-authored more than 100 articles in refereed journals and conference proceedings. Dr. Rontogiannis has been a program committee member in more than 25 conferences, in one of them as Co-Chair and in three of them as Area Chair. He has served at the Editorial Boards of the EURASIP Journal on Advances in Signal Processing, Springer (2008–2017), and the EURASIP Signal Processing Journal, Elsevier (since 2011). Since 2017, he has been serving as an Associate Editor of the IEEE Transactions on Signal Processing. He is a member of the EURASIP and the Technical Chamber of Greece. Since January 2018, he has also been the Chair of the IEEE Signal Processing Society Greece Chapter.
\end{IEEEbiography}

\begin{IEEEbiography}{Eleftherios Kofidis}
received the Diploma and Ph.D. degrees both from the Department of Computer Engineering and Informatics, University of Patras, Patras, Greece, in 1990 and 1996, respectively. From 1996 to 1998 he served in the Hellenic Army. From 1998 to 2000, he was a Postdoctoral Fellow at the Institut National des T\'{e}l\'{e}communications (INT), \'{E}vry, France (now T\'{e}l\'{e}com SudParis). From 2001 to 2004, he was a Research Associate at the University of Athens, and Adjunct Professor at the Universities of Peloponnese and Piraeus, Greece. In 2004, he joined the Dept. of Statistics and Insurance Science, University of Piraeus, Greece, where he is now Associate Professor. He is also affiliated with the Computer Technology Institute \& Press “Diophantus” (CTI), Greece. His research interests are in signal processing and machine learning, with applications including communications and medical imaging. Dr. Kofidis has served as technical program co-chair in two international conferences (CIP-2008 and DSP-2009) and as a technical program committee member and reviewer in a number of conferences and journals. He has served as Associate Editor in the IEEE Transactions on Signal Processing, the EURASIP Journal on Advances in Signal Processing (JASP) and the IET Signal Processing journal. He co-organized three special sessions on filter bank-based multicarrier systems (ISWCS-2012, EW-2014, SPAWC-2015) and was Lead Guest Editor for a JASP special issue on this subject.
\end{IEEEbiography}

\begin{IEEEbiography}{Paris V. Giampouras}
was born in Athens, Greece, in 1986. In 2011, he received the Diploma degree in electrical and computer engineering from the National Technical University of Athens, Athens, Greece, and the M.Sc. degree in information technologies in medicine and biology in 2014 from the Department of Informatics and Telecommunications, National and Kapodistrian University of Athens, Athens, Greece, where he received the Ph.D. degree in 2018. Since 2019 he has been a postdoctoral fellow at the  Mathematical Institute for Data Science, Johns Hopkins University, Baltimore, MD. His main research interests include the areas of signal processing and machine learning focusing on sparse and low-rank representations of large-scale data and their application to collaborative filtering and hyper-spectral image processing.
\end{IEEEbiography}

\end{document}